# Reproducing kernel Hilbert spaces of Gaussian priors

## A. W. van der Vaart[1] and J. H. van Zanten[*,1]

*Vrije Universiteit Amsterdam*

**Abstract:** We review definitions and properties of reproducing kernel Hilbert spaces attached to Gaussian variables and processes, with a view to applications in nonparametric Bayesian statistics using Gaussian priors. The rate of contraction of posterior distributions based on Gaussian priors can be described through a concentration function that is expressed in the reproducing Hilbert space. Absolute continuity of Gaussian measures and concentration inequalities play an important role in understanding and deriving this result. Series expansions of Gaussian variables and transformations of their reproducing kernel Hilbert spaces under linear maps are useful tools to compute the concentration function.

## Contents

## 1. Introduction

Ghosal, Ghosh and van der Vaart considered in [4] the rate of contraction of a posterior distribution based on i.i.d. observations to the true density. Given prior probability measures $\Pi_n$ defined on a set $\mathcal{P}$ of densities $p$ relative to a given $\sigma$-finite measure on a measurable space (such that the maps $(x,p) \mapsto p(x)$ are jointly

*Supported in part by the Netherlands Organization for Scientific Research NWO.
[1] Department of Mathematics, Vrije Universiteit, De Boelelaan 1081a, 1081 HV Amsterdam, The Netherlands, e-mail: aad@cs.vu.nl; harry@cs.vu.nl

*AMS 2000 subject classifications:* 60G15, 62G05.
*Keywords and phrases:* Bayesian inference, rate of convergence.

measurable) and observations $X_1, \ldots, X_n$, they characterized the rate $\varepsilon_n \downarrow 0$ at which the posterior distribution

$$\Pi_n(B | X_1, \ldots, X_n) = \frac{\int_B \prod_{i=1}^n p(X_i) \, d\Pi_n(p)}{\int \prod_{i=1}^n p(X_i) \, d\Pi_n(p)} \tag{1.1}$$

contracts to $p_0$ if the observations are an i.i.d. sample from this density, i.e. the rate for which

$$E_{p_0} \Pi_n \big( p: d(p, p_0) > M\varepsilon_n | X_1, \ldots, X_n \big) \to 0,$$

for sufficiently large $M$. In their results $d$ can be the Hellinger distance, the $L_1$-distance, or the $L_2$-distance if the densities are uniformly bounded above.

The paper [15] applied these results to priors $\Pi_n$ constructed from Gaussian processes. They consider a prior $\Pi_n$ constructed as the distribution of $p_W$, for $W$ a Gaussian random element in a Banach space $(\mathbb{B}, \|\cdot\|)$ and $w \mapsto p_w$ a map such that, for some constant $C$ and all $v, w \in \mathbb{B}$ with $\|v - w\|$ bounded above by some fixed constant,

$$\begin{aligned} d(p_v, p_w) &\leq C\|v - w\|, \\ K(p_v, p_w) &\leq C\|v - w\|^2, \\ V(p_v, p_w) &\leq C\|v - w\|^2. \end{aligned}$$

Here $K(p, q) = \int \log(p/q) \, p \, d\mu$ is the Kullback–Leibler divergence and $V(p, q) = \int \big(\log(p/q)\big)^2 p \, d\mu$. This setting covers, for instance, the case of density estimation on $[0, 1]$ as considered in Tokdar and Ghosh [14], with $d$ the Hellinger distance, the Banach space equal to $\mathbb{B} = C[0, 1]$ and

$$p_w(x) = \frac{e^{w_x}}{\int_0^1 e^{w_y} \, dy}.$$

It also covers logistic or probit regression as considered in [5] with appropriate choices and several other situations, as shown in [15].

In the latter paper it is shown that if the true density takes the form $p_0 = p_{w_0}$, then the rate of posterior contraction $\varepsilon_n$ is characterized by the pair of equations

$$\inf_{h \in \mathbb{H}: \|h - w_0\| < \varepsilon_n} \|h\|_{\mathbb{H}}^2 \leq n\varepsilon_n^2, \tag{1.2}$$

$$-\log P\big(\|W\| < \varepsilon_n\big) \leq n\varepsilon_n^2. \tag{1.3}$$

Here $(\mathbb{H}, \|\cdot\|_{\mathbb{H}})$ is the *reproducing kernel Hilbert space* (RKHS) of the Gaussian variable, and $P(\|W\| < \varepsilon)$ is its *small ball probability* (cf. [11]). Both equations have a minimal solution $\varepsilon_n$, and the rate is the worse of the two solutions. The second depends only on the prior, and gives a maximal rate regardless of the true parameter $w_0$, whereas the first involves the true parameter.

The reproducing kernel Hilbert space arises because it determines the support and the "geometry" of the concentration of the Gaussian measure, which are crucial for its success as a prior. Results on RKHSs of Gaussian variables are spread over many research papers, and sometimes seem to belong to what is "well known" without clear references. Moreover, there are different definitions for stochastic processes and Borel measurable maps in a separable Banach space. In this paper we review definitions, investigate when the different definitions agree, and derive results that are useful for the construction of priors and the study of posterior distributions.

## 2. Definitions and elementary properties

In this section we give and compare two definitions of RKHS, one for stochastic processes and one for Borel measurable maps in a Banach space.

### 2.1. Gaussian processes

A zero-mean Gaussian stochastic process $W = (W_t : t \in T)$ is a set of random variables $W_t$ indexed by an arbitrary set $T$ and defined on a common probability space $(\Omega, \mathcal{U}, P)$ such that each finite subset possesses a zero-mean multivariate normal distribution. The finite-dimensional distributions of such a process are determined by the covariance function $K : T \times T \to \mathcal{R}$, defined by

$$K(s,t) = EW_s W_t.$$

The *reproducing kernel Hilbert space (RKHS)* attached to the Gaussian process $W$ is the completion $\mathbb{H}$ of the linear space of all functions

(2.1) $\quad t \mapsto \sum_{i=1}^{k} \alpha_i K(s_i, t), \qquad \alpha_1, \ldots, \alpha_k \in \mathcal{R}, s_1, \ldots, s_k \in T, k \in \mathbb{N},$

relative to the norm induced by the inner product

(2.2) $\quad \left\langle \sum_{i=1}^{k} \alpha_i K(s_i, \cdot), \sum_{j=1}^{l} \beta_j K(t_j, \cdot) \right\rangle_{\mathbb{H}} = \sum_{i=1}^{k} \sum_{j=1}^{l} \alpha_i \beta_j K(s_i, t_j).$

It can be checked that this definition is independent of the representation of the functions on the left, and that this defines a valid inner product.

The completion of the collection of functions (2.1) is an abstract metric-topological operation using the metric induced by the inner product (2.2) only. As such the completion is not a space of functions $f : T \to \mathcal{R}$. However, it can be identified with a space of functions $f : T \mapsto \mathcal{R}$, through the *reproducing formula*

$$f(t) = \langle f, K(t, \cdot) \rangle_{\mathbb{H}}.$$

For $f$ a linear combination of the form $\sum_{i=1}^{k} \alpha_i K(s_i, \cdot)$ this formula follows from the definition (2.2) of the inner product $\langle \cdot, \cdot \rangle_{\mathbb{H}}$. For general $f \in \mathbb{H}$ the (extended) inner product on the right (with the extended function $K(t, \cdot)$) is well defined through the completion operation, and can be used to define a function $f : T \mapsto \mathcal{R}$.

Alternatively, the function in (2.1) can be written as

(2.3) $\quad t \mapsto EW_t H, \qquad H = \sum_i \alpha_i W_{s_i}.$

With the function in the display written as $EW.H$, the inner product (2.2) is equal to

$$\langle EW.H_1, EW.H_2 \rangle_{\mathbb{H}} = EH_1 H_2.$$

Thus the map $H \mapsto EW.H$ is an isometry for the norm of the $L_2$-space attached to the probability space $(\Omega, \mathcal{U}, P)$ on which the process $W$ is defined and the RKHS-norm. The stochastic process RKHS $\mathbb{H}$, which is defined as the completion of the set of functions (2.3), is therefore precisely the set of functions $t \mapsto EW_t H$ with $H$ ranging over the closure of the set of linear combinations $H = \sum_i \alpha_i W_{s_i}$ in $L_2(\Omega, \mathcal{U}, P)$ (known as the *first order chaos of $W$*). It follows again that we can view $\mathbb{H}$ as a Hilbert space of functions on $T$.

### 2.2. Gaussian elements in a Banach space

A Borel measurable random element $W$ with values in a separable Banach space $(\mathbb{B}, \|\cdot\|)$ is called *Gaussian* if the random variable $b^*W$ is normally distributed for any element $b^*$ of the dual space $\mathbb{B}^*$ of $\mathbb{B}$, and it is called *zero-mean* if the mean of every such variable $b^*W$ is zero. Henceforth we shall only consider zero-mean Gaussian variables.

It is well known that the norm $\|W\|$ of a zero-mean Gaussian variable, which is a finite random variable by the assumption that $W$ takes its values in $\mathbb{B}$, has sub-Gaussian tails. (cf. Corollary 5.1 below, or, e.g., [17], Propositions A.2.1 and A.2.3, for a direct proof.) In particular, all moments $\mathrm{E}\|W\|^p$ are finite. We set

$$\sigma^2(W) = \sup_{b^* \in \mathbb{B}^*: \|b^*\|=1} \mathrm{E}b^*(W)^2.$$

This is a finite number, bounded by $\mathrm{E}\|W\|^2$.

For every element $b^* \in \mathbb{B}^*$ we define $Sb^* \in \mathbb{B}$ as the Pettis integral $\mathrm{E}Wb^*(W)$ of the $\mathbb{B}$-valued random element $Wb^*(W)$. By definition, this *Pettis integral* is an element $Sb^*$ of $\mathbb{B}$ such that $b_2^*(Sb^*) = \mathrm{E}b_2^*(W)b^*(W)$ for every $b_2^* \in \mathbb{B}^*$. The following lemma allows us to derive the existence of the Pettis integral from the fact that $\mathrm{E}\|W\|^2 < \infty$.

**Lemma 2.1.** *If $X$ is a Borel measurable map in a separable Banach space $\mathbb{B}$ with $\mathrm{E}\|X\| < \infty$, then there exists an element $b \in \mathbb{B}$ such that $b^*(b) = \mathrm{E}b^*(X)$ for every $b^* \in \mathbb{B}^*$.*

*Proof.* Because the Banach space is assumed separable, the map $X$ is automatically tight (e.g. [17], 1.3.2). Therefore, for any $n \in \mathbb{N}$ there exists a compact set $K$ such that $\mathrm{E}\|X\|1_{X \notin K} < 1/n$. This compact set can be partitioned into finitely many sets $B_i$ of diameter smaller than $1/n$. Without loss of generality these partitions can be chosen as successive refinements for increasing $n$. Let $X_n = \sum_i b_i 1_{X \in B_i}$ for $b_i$ arbitrary points in the partitioning sets. Then $\mathrm{E}X_n := \sum_i b_i \mathrm{P}(X \in B_i)$ satisfies $b^*(\mathrm{E}X_n) = \mathrm{E}b^*(X_n)$ for every $b^* \in \mathbb{B}^*$. Furthermore, the sequence $\mathrm{E}X_n$ is a Cauchy sequence in $\mathbb{B}$, because $\|\mathrm{E}X_n - \mathrm{E}X_m\| = \sup_{\|b^*\|=1} |\mathrm{E}b^*(X_n - X_m)| \le \mathrm{E}\|X_n - X_m\| \to 0$ as $n, m \to \infty$. Because $\mathrm{E}\|X_n - X\| < 2/n$, we have that $b^*(\mathrm{E}X_n) = \mathrm{E}b^*(X_n) \to \mathrm{E}b^*(X)$ for every $b^* \in \mathbb{B}$. The strong limit $b$ of the sequence $\mathrm{E}X_n$ is of course also a weak limit, whence $b^*(b) = \mathrm{E}b^*(X)$ for every $b^* \in \mathbb{B}^*$. □

The *reproducing kernel Hilbert space (RKHS)* $\mathbb{H}$ attached to $W$ is the completion of the range $S\mathbb{B}^*$ of the map $S: \mathbb{B}^* \to \mathbb{B}$ defined by $Sb^* = \mathrm{E}Wb^*(W)$ for the inner product

$$\langle Sb_1^*, Sb_2^* \rangle_\mathbb{H} = \mathrm{E}b_1^*(W)b_2^*(W).$$

By the Hahn–Banach theorem and the Cauchy–Schwarz inequality,

$$\|Sb^*\| = \sup_{b_2^* \in \mathbb{B}^*: \|b_2^*\|=1} |b_2^*(Sb^*)| = \sup_{b_2^* \in \mathbb{B}^*: \|b_2^*\|=1} |\mathrm{E}b_2^*(W)b^*(W)|$$

(2.4) $$\le \sigma(W)\bigl(\mathrm{E}b^*(W)^2\bigr)^{1/2} = \sigma(W)\|Sb^*\|_\mathbb{H}.$$

It follows that the RKHS-norm on the set $S\mathbb{B}^*$ is stronger than the original norm, so that a $\|\cdot\|_\mathbb{H}$-Cauchy sequence in $S\mathbb{B}^* \subset \mathbb{B}$ is a $\|\cdot\|$-Cauchy sequence in $\mathbb{B}$. Consequently, the RKHS, which is by definition the completion of the set $S\mathbb{B}^*$

under the RKHS norm, can be identified with a subset of $\mathbb{B}$. In terms of the unit balls $\mathbb{B}_1$ and $\mathbb{H}_1$ of $\mathbb{B}$ and $\mathbb{H}$ the preceding display can be written as

$$\mathbb{H}_1 \subset \sigma(W)\mathbb{B}_1. \tag{2.5}$$

In other words, the norm of the embedding $i\colon \mathbb{H} \to \mathbb{B}$ is bounded by $\sigma(W)$.

**Lemma 2.2.** *The map* $S\colon \mathbb{B}^* \to \mathbb{H}$ *is weak-\* continuous.*

*Proof.* The unit ball $\mathbb{B}_1^*$ of the dual space is weak-\* metrizable ([12], 3.16). Therefore the restricted map $S\colon \mathbb{B}_1^* \to \mathbb{H}$ is weak-\* continuous if and only if weak-\* convergence of a sequence $b_n^*$ in $\mathbb{B}_1^*$ to an element $b^*$ implies that $Sb_n^* \to Sb^*$ in $\mathbb{H}$. Now the weak-\* convergence $b_n^* \to b^*$ is by definition pointwise convergence on $\mathbb{B}$. Then the sequence $(b_n^* - b^*)(W)$ tends to zero (almost) surely, and hence also in distribution. Because each of these variables is zero-mean Gaussian, this implies that the variances tend to zero, i.e. $\|Sb_n^* - Sb^*\|_{\mathbb{H}}^2 = \mathrm{E}(b_n^* - b)^2(W) \to 0$. (Alternatively, use the uniform integrability of the variables $b^*W$ instead of the Gaussianity.)

This concludes the proof that the restriction of $S$ to the unit ball $\mathbb{B}_1^*$ is continuous. A weak-\* converging net $b_n^*$ in $\mathbb{B}^*$ is necessarily bounded in norm, by the Banach–Steinhaus theorem ([12], 2.5), and hence is contained in a multiple of the unit ball. The continuity of the restriction then shows that $Sb_n^* \to Sb^*$, which concludes the proof. □

**Corollary 2.1.** *If* $\mathbb{B}_0^*$ *is a weak-\* dense subset of* $\mathbb{B}^*$, *then* $\mathbb{H}$ *is the completion of* $S\mathbb{B}_0^*$.

By the definitions $\langle Sb^*, S\underline{b}^* \rangle_{\mathbb{H}} = \mathrm{E}b^*W\underline{b}^*W = b^*(S\underline{b}^*)$, for any $b^*, \underline{b}^* \in \mathbb{B}^*$. By continuity of the inner product this extends to the *reproducing formula*:

$$\langle Sb^*, h \rangle_{\mathbb{H}} = b^*(h), \tag{2.6}$$

which is valid for every $h \in \mathbb{H}$ and $b^* \in \mathbb{B}^*$.

Just as for stochastic processes there is an alternative representation of the RKHS through "first chaos", in the present setting defined as the closed linear span of the variables $b^*W$ in $L_2(\Omega, \mathcal{U}, \mathrm{P})$. The elements $Sb^*$ of the RKHS can be written $Sb^* = \mathrm{E}HW$ for $H = b^*W$, and the RKHS-norm of $Sb^*$ is by definition the $L_2(\Omega, \mathcal{U}, \mathrm{P})$-norm of this $H$. This immediately implies the following lemma. Note that $\mathrm{E}HW$ is well defined as a Pettis integral for every $H \in L_2(\Omega, \mathcal{U}, \mathrm{P})$, by Lemma 2.1.

**Lemma 2.3.** *The RKHS is the set of Pettis integrals* $\mathrm{E}HW$ *for $H$ ranging over the closed linear span of the variables $b^*W$ in $L_2(\Omega, \mathcal{U}, \mathrm{P})$ with inner product* $\langle \mathrm{E}H_1W, H_2W \rangle_{\mathbb{H}} = \mathrm{E}H_1H_2$.

It is useful to decompose the map $S\colon \mathbb{B}^* \to \mathbb{B}$ as $S = A^*A$ for $A^*\colon L_2(\Omega, \mathcal{U}, \mathrm{P}) \to \mathbb{B}$ and $A\colon \mathbb{B}^* \to L_2(\Omega, \mathcal{U}, \mathrm{P})$ given by

$$\begin{aligned} A^*H &= \mathrm{E}HW, \\ Ab^* &= b^*W. \end{aligned}$$

It may be checked that the operators $A$ and $A^*$ are indeed adjoints, after identifying $\mathbb{B}$ with a subset of its second dual space $\mathbb{B}^{**}$ under the canonical embedding ([12], 3.15, 4.5), as the notation suggests. By the preceding lemma the RKHS is the image of the first chaos space under $A^*$. Because $R(A)^\perp = N(A^*)$ the full range $R(A^*) = A^*\big(L_2(\Omega, \mathcal{U}, \mathrm{P})\big)$ is not bigger than the image of the first chaos, although the map $A^*: L_2(\Omega, \mathcal{U}, P) \to \mathbb{H}$ is an isometry only if restricted to the first chaos space.

Recall that an operator is compact if it maps bounded sets into precompact sets, or, equivalently, maps bounded sequences into sequences that possess a converging subsequence.

**Lemma 2.4.** *The maps $A^*: L_2(\Omega, \mathcal{U}, P) \to \mathbb{B}$ and $A: \mathbb{B}^* \to L_2(\Omega, \mathcal{U}, P)$ and $S: \mathbb{B}^* \to \mathbb{B}$ are compact for the norms.*

*Proof.* In general an operator is compact if and only if its adjoint is compact, and a composition with a compact operator is compact (see ([12], 4.19). To prove the compactness of $A$ fix some sequence $b_n^*$ in the unit ball $\mathbb{B}_1^*$. As the unit ball is weak-* compact by the Banach–Alaoglu theorem ([12], 4.3(c)), there exists a subsequence along which $b_{n_j}^*$ converges pointwise on $\mathbb{B}$ to a limit $b^*$. Consequently $b_{n_j}^*(W) \to b^*(W)$ almost surely, and hence in second mean. □

As a consequence we can conclude that the unit ball of the RKHS is precompact in $\mathbb{B}$. Indeed, $\mathbb{H}_1 = A^* \mathbb{U}_1$ for $\mathbb{U}_1$ the unit ball of $L_2(\Omega, \mathcal{U}, P)$, and hence is precompact by the compactness of $A^*$.

**Example 2.1** (Hilbert space). The *covariance operator* of a mean zero Gaussian random element $W$ in a Hilbert space $\mathbb{B}$ with inner product $\langle \cdot, \cdot \rangle$ is the map $S: \mathbb{B} \to \mathbb{B}$ that satisfies $E\langle W, b_1 \rangle \langle W, b_2 \rangle = \langle b_1, S b_2 \rangle$. It is well known that $S$ is continuous, linear, positive, self-adjoint, and of finite trace, and hence it possesses a square root, which is another positive, self-adjoint operator $S^{1/2}: \mathbb{B} \to \mathbb{B}$ such that $S^{1/2} S^{1/2} = S$. (The square root can also be described as having the same eigenfunctions as $S$ with eigenvalues the square roots of the eigenvalues of $S$.) The RKHS of $W$ can be characterized as the range of $S^{1/2}$ equipped with the norm $\|S^{1/2} b\|_{\mathbb{H}} = \|b\|$.

To see this note that the covariance operator $S$ is exactly the operator $S$ as defined previously, after the usual identification of the dual space $\mathbb{B}^*$ with $\mathbb{B}$ itself: $b \in \mathbb{B}$ corresponds to the element $b_1 \mapsto \langle b, b_1 \rangle$ of $\mathbb{B}^*$. Hence the RKHS is the completion of the elements $Sb$ under the square norm $\|Sb\|_{\mathbb{H}}^2 = E\langle W, b \rangle^2 = \langle b, Sb \rangle = \|S^{1/2} b\|^2$. This is the same as the completion of the set of functions $S^{1/2} c$ (with $c = S^{1/2} b$) under the norm $\|S^{1/2} c\|_V^2 = \|c\|^2$. The latter set is of course already complete, so that completion is superfluous.

## 2.3. Comparison

If the sample paths $t \mapsto W_t$ of a stochastic process $W = (W_t : t \in T)$ belong to a Banach space of functions, then the process can be viewed as a map $W$ into the Banach space. If it is a Borel measurable map, then the preceding gives two definitions of a RKHS. The two definitions will coincide provided the dual space can be appropriately related to the covariance function. In particular, if the coordinate projections $\pi_t: \mathbb{B} \to \mathcal{R}$, defined by $b \mapsto b(t)$, are elements of the dual space, then $W_t = \pi_t(W)$ and the covariance function $K(s, t) = E W_s W_t$ takes the form $E \pi_s(W) \pi_t(W) = \langle S \pi_s, S \pi_t \rangle_{\mathbb{H}}$. If the other elements $S b^*$ are determined by the elements $S \pi_t$, then the two definitions should be the same. It appears that in general some conditions are needed to make the link between the two definitions. For the Banach space $\ell^\infty(T)$ of uniformly bounded functions $z: T \to \mathcal{R}$ equipped with the uniform norm $\|z\| = \sup\{|z(t)| : t \in T\}$, this can always be done.

The following result is probably known to the experts, but we do not know a published reference.

**Theorem 2.1.** *If $W$ is a Borel measurable zero-mean Gaussian random element in a complete separable subspace of $\ell^\infty(T)$ equipped with the uniform norm, then the Banach space RKHS and the stochastic process RKHS coincide. Furthermore $S\pi_t = K(t, \cdot)$.*

*Proof.* For a given tight Borel measurable random element $W$ in $\ell^\infty(T)$ there exists a semimetric $\rho$ on $T$ under which $T$ is totally bounded and such that $W$ takes its values in the subspace $UC(T, \rho)$ of functions $f: T \to \mathcal{R}$ that are uniformly continuous relative to $\rho$ (e.g. [17], Lemma 1.5.9). Thus we may assume without loss of generality that $W$ takes its values in $UC(T, \rho)$ for such a semimetric $\rho$. The space $UC(T, \rho)$ is a Banach space under the supremum norm $\|f\| = \sup\{|f(t)|: t \in T\}$. Let $K(s, t) = \mathrm{E}W_s W_t$.

The coordinate projections $\pi_t: f \mapsto f(t)$ belong to the dual space $UC(T, \rho)^*$. The corresponding Pettis integral $S\pi_t$ is the function $K(t, \cdot)$. This follows because it is contained in $UC(T, \rho)$ and, furthermore, for every $s \in T$,

$$\pi_s\bigl(K(t, \cdot)\bigr) = K(t, s) = \mathrm{E}W_s W_t = \mathrm{E}\pi_s(W)\pi_t(W).$$

Because the coordinate projections $\pi_t f$ identify $f$ uniquely it follows that $K(t, \cdot) = \mathrm{E}W\pi_t(W) = S\pi_t$.

Thus the stochastic process RKHS, defined as the completion of the linear combinations (2.1), is contained in the Banach space RKHS. The inner products on the two spaces agree, because

$$\langle S\pi_s, S\pi_t \rangle_\mathbb{H} = \mathrm{E}\pi_t(W)\pi_s(W) = K(s, t) = \bigl\langle K(s, \cdot), K(t, \cdot) \bigr\rangle_\mathbb{H}.$$

By the Riesz representation theorem an arbitary element of $UC(T, \rho)^*$ is a map $f \mapsto \int \bar{f}(t)\,d\mu(t)$ for a signed Borel measure on the completion $\bar{T}$ of $T$ and $\bar{f}: \bar{T} \to \mathcal{R}$ is the continuous extension of $f$. Because $T$ is totally bounded we can write it for each $m \in \mathbb{N}$ as a finite union of sets of diameter smaller than $1/m$. If we define $\mu_m$ as the measure obtained by concentrating the masses of $\mu$ on the partitioning sets in a fixed, single point in the partitioning set, then $\int \bar{f}\,d\mu_m \to \int \bar{f}\,d\mu$ as $m \to \infty$, for each $f \in UC(T, \rho)$. The map $f \mapsto \int \bar{f}\,d\mu_m$ is a linear combination of coordinate projections. It follows that for any $b^* \in UC(T, \rho)^*$ there exists a sequence $b_m^*$ of linear combinations of coordinate projections that converges pointwise on $UC(T, \rho)$ to $b^*$. In other words, the linear span $\mathbb{B}_0^*$ of coordinate projections is weak-* dense in $UC(T, \rho)^*$, and hence the RKHS is the completion of $S\mathbb{B}_0^*$, by Lemma 2.1. □

**Example 2.2.** The preceding theorem applies, for instance, to the space of continuous functions $z: T \to \mathcal{R}$ on a compact metric space $T$. For instance $C[0, 1]$.

A more general connection between the two definitions of a RKHS can be made by embedding the Banach space $\mathbb{B}$ in its second dual (see [12], 4.15). This is somewhat technical and will not be needed in the rest of the paper. The canonical embedding is, as usual, the identification of $b \in \mathbb{B}$ with the map $b^{**}: \mathbb{B}^* \to \mathcal{R}$ defined by $b^{**}(b^*) = b^*(b)$. A Borel measurable random element $W$ in $\mathbb{B}$ becomes identified in this way with the stochastic process $W^{**} = \bigl(b^*(W): b^* \in \mathbb{B}^*\bigr)$, which has covariance function

$$K(b_1^*, b_2^*) = \mathrm{E}b_1^*(W)b_2^*(W).$$

The stochastic process RKHS $\mathbb{H}$ attached to this process in Section 2.1 is the completion of the set of functions $K(b^*, \cdot): \mathbb{B}^* \to \mathcal{R}$ relative to the inner product

$$\langle K(b_1^*, \cdot), K(b_2^*, \cdot)\rangle_\mathbb{H} = K(b_1^*, b_2^*) = \mathrm{E}b_1^*(W)b_2^*(W).$$

The function $K(b^*, \cdot)$ is exacly the Pettis integral $\mathrm{E}Wb^*(W)$, written $Sb^*$ in the preceding and now viewed as an element of $\mathbb{B}^{**}$; and the inner product in the display is exactly $\langle Sb_1^*, Sb_2^* \rangle_\mathbb{H}$. Thus the two definitions of RKHS coincide, after identification of $\mathbb{B}$ and its image in $\mathbb{B}^{**}$ under the canonical embedding.

## 3. Absolute continuity

Given a zero-mean Gaussian process $W = (W_t : t \in T)$ with covariance kernel $K$ defined on a probability space $(\Omega, \mathcal{U}, \mathrm{P})$ with RKHS $\mathbb{H}$ as defined in Section 2.1, we can define a map $U : \mathbb{H} \to L_2(\Omega, \mathcal{U}, \mathrm{P})$ by defining

$$UK(t, \cdot) = W_t, \tag{3.1}$$

and extending linearly and continuously. This map is an Hilbert space isometry, since
$$\mathrm{E}UK(s, \cdot)UK(t, \cdot) = \mathrm{E}W_s W_t = K(s,t) = \langle K(s, \cdot), K(t, \cdot)\rangle_\mathbb{H}.$$
This isometry property also implies the existence of the extension. It follows that the process $(Uh : h \in \mathbb{H})$ is the *iso-Gaussian process* indexed by $\mathbb{H}$: a mean-zero Gaussian process with covariance function $\mathrm{E}UgUh = \langle g, h \rangle_\mathbb{H}$.

The process $W$ induces a distribution $P^W$ on the product $\sigma$-field of $\mathcal{R}^T$. For a function $f : T \mapsto \mathcal{R}$ the process $(W_t + f(t) : t \in T)$ induces another distribution $P^{W+f}$ on the same space.

**Lemma 3.1.** *If $f \in \mathbb{H}$, then $P^{W+f}$ and $P^W$ are equivalent and*

$$\frac{dP^{W+f}}{dP^W}(W) = e^{Uf - \frac{1}{2}\|f\|_\mathbb{H}^2}, \qquad \text{a.s.}$$

*Proof.* The process $W$ is the "subprocess" $W^\mathbb{G} = (Ug : g \in \mathbb{G})$ of the iso-Gaussian process $W^\mathbb{H} = (Uh : h \in \mathbb{H})$ for $\mathbb{G}$ the set of functions $K(t, \cdot)$ with $t$ ranging over $T$. From the general theory of Gaussian processes

$$\frac{dP^{W^\mathbb{H} + (\langle h, f\rangle_\mathbb{H} : h \in \mathbb{H})}}{dP^{W^\mathbb{H}}}(W^\mathbb{H}) = e^{Uf - \frac{1}{2}\|f\|_\mathbb{H}^2}, \qquad \text{a.s.} \tag{3.2}$$

The process $W^\mathbb{G}$ arises from the iso-Gaussian process by the projection $\pi_\mathbb{G} : \mathcal{R}^\mathbb{H} \to \mathcal{R}^\mathbb{G}$. The corresponding Radon–Nikodym derivative can be found as the conditional expectation

$$\frac{dP^{W^\mathbb{G} + (\langle g, f\rangle_\mathbb{H} : g \in \mathbb{G})}}{dP^{W^\mathbb{G}}}(W^\mathbb{G}) = \mathrm{E}\Big(\frac{dP^{W^\mathbb{H} + (\langle h, f\rangle_\mathbb{H} : h \in \mathbb{H})}}{dP^{W^\mathbb{H}}}(W^\mathbb{H}) \,\big|\, W^\mathbb{G}\Big).$$

Because $\mathrm{lin}\,(\mathbb{G})$ is dense in $\mathbb{H}$ by construction and $U$ is continuous, the variable $Uf$ is the $L_2(\Omega, \mathcal{U}, \mathrm{P})$-limit of a sequence $Ug_n$ with $(g_n) \subset \mathrm{lin}\,(\mathbb{G})$ and hence is measurable relative to the completion of the $\sigma$-field generated by $W^\mathbb{G}$. Consequently, the right side of (3.2) is $\mathcal{W}^\mathbb{G}$-measurable as well and hence the conditional expectation in the preceding display is unnecessary.

Finally, note that the shift $\langle g, f \rangle_\mathbb{H}$ is exactly the function $f$ after the identification $g \leftrightarrow K(t, \cdot)$, by the reproducing property: $f(t) = \langle K(t, \cdot), f \rangle_\mathbb{H}$ for every $t \in T$. □

Let $\mathbb{H}$ be the abstract RKHS attached to a zero-mean, Borel measurable, Gaussian random element $W$ in a separable Banach space $\mathbb{B}$ defined on a probability space $(\Omega, \mathcal{U}, \mathrm{P})$. Let $U : \mathbb{H} \to L_2(\Omega, \mathcal{U}, \mathrm{P})$ be the isometry defined by

$$U(Sb^*) = b^*(W), \qquad b^* \in \mathbb{B}^*, \tag{3.3}$$

and extending continuously. It is the same map $U$ as in (3.1) if we make the identification $S\pi_t = K(t, \cdot)$ of Theorem 2.1; also $US = A$ for $A$ defined in Section 2.2. As before the map $U$ is an isometry. The preceding lemma can be translated to the present situation.

**Lemma 3.2.** *If $h \in \mathbb{H}$ then the distributions $P^{W+h}$ and $P^W$ of $W + h$ and $W$ on $\mathbb{B}$ are equivalent and*

$$\frac{dP^{W+h}}{dP^W}(W) = e^{Uh - \frac{1}{2}\|h\|_{\mathbb{H}}^2}, \qquad \text{a.s.}$$

*Proof.* The process $W^{**} = \big(b^*(W)\colon b^* \in \mathbb{B}^*\big)$ arising from $W$ through the canonical embedding generates the same $\sigma$-field on the underlying probability space as $W$ and can be viewed as a measurable transformation of $W$ under the map $\phi\colon \mathbb{B} \to \mathcal{R}^{\mathbb{B}^*}$ given by $\phi(b)(b^*) = b^*(b)$. The process $W + h$ is transformed in the process $W^{**} + h^{**} = \phi(W + h)$. The result therefore follows from Lemma 3.1.

The following alternative proof is given in Proposition 2.1 in [3]. The isometry property of $U$ shows that $\mathrm{E}(Uh)^2 = \|h\|_{\mathbb{H}}^2$. Because $Uh$ is in the closed linear span of the zero-mean Gaussian variables $USb^* = b^*W$, it is itself zero-mean Gaussian. It follows that

$$dQ = e^{Uh - \frac{1}{2}\|h\|_{\mathbb{H}}^2} dP$$

defines a probability measure on $(\Omega, \mathcal{U})$. For any $b_1^*, b_2^* \in \mathbb{B}^*$ the joint distribution of $(USb_1^*, USb_2^*) = (b_1^*W, b_2^*W)$ is bivariate normal with mean zero and covariance matrix $\big(\langle Sb_i^*, Sb_j^* \rangle_{\mathbb{H}}\big)_{i,j=1,2}$. By taking limits we see that for every $h \in \mathbb{H}$ the joint distribution of $(b_1^*W, Uh)$ is bivariate normal with mean zero and covariance matrix $\Sigma$ with $\Sigma_{1,1} = \|Sb_1^*\|_{\mathbb{H}}^2$, $\Sigma_{1,2} = \langle Sb_1^*, h \rangle_{\mathbb{H}}$ and $\Sigma_{2,2} = \|h\|_{\mathbb{H}}^2$. Thus

$$\mathrm{E}_Q e^{ib_1^*W} = \mathrm{E}e^{ib_1^*W} e^{Uh - \frac{1}{2}\|h\|_{\mathbb{H}}^2} = e^{\frac{1}{2}(i,1)\Sigma(i,1)^T} e^{-\frac{1}{2}\|h\|_{\mathbb{H}}^2} = e^{-\frac{1}{2}\Sigma_{1,1} + i\Sigma_{1,2}}.$$

The right side is also equal to

$$\mathrm{E}e^{ib_1^*W + i\langle Sb_1^*, h\rangle_{\mathbb{H}}} = \mathrm{E}e^{ib_1^*(W+h)}.$$

The last step follows from the reproducing formula (2.6). We conclude that the distribution of $W + h$ under P is the same as the distribution of $W$ under $Q$, i.e. $\mathrm{P}(W + h \in B) = \mathrm{E}_Q 1_B(W) = \mathrm{E}1_B(W)(dQ/dP)$. □

The preceding lemma requires that the shift $h$ is contained in the RKHS. If this is not the case, then there is no density.

**Lemma 3.3.** *If $b \notin \mathbb{H}$ then the distributions $P^{W+b}$ and $P^W$ of $W + b$ and $W$ on $\mathbb{B}$ are orthogonal.*

*Proof.* By Lemma 5.1 (below) the closure $\bar{\mathbb{H}}$ of $\mathbb{H}$ in $\mathbb{B}$ is the support of $W$. Because the affine spaces $\bar{\mathbb{H}}$ and $\bar{\mathbb{H}} + b$ are disjoint if $b \notin \bar{\mathbb{H}}$, the assertion is clear if $b \in \mathbb{B} - \bar{\mathbb{H}}$. Therefore, it is not a loss of generality to assume that $\mathbb{B}$ is the closure of $\mathbb{H}$.

Fix a sequence $\{b_n^*\} \subset \mathbb{B}^*$ whose linear span is dense (for the norm) in $\mathbb{B}^*$ and is such that the variables $b_n^*W$ are i.i.d. standard normal variables. We prove the existence of such a sequence at the end of the proof. We claim that $\mathbb{H} = \{b \in \mathbb{B}: \sum_{n=1}^\infty (b_n^*b)^2 < \infty\}$. Indeed, the sequence $h_n = Sb_n^*$ is orthonormal in $\mathbb{H}$ by the definition of the inner product in $\mathbb{H}$ and $\mathrm{lin}\,(h_n) = S\,\mathrm{lin}\,(b_n^*)$ is dense in $S\mathbb{B}^*$ by construction of the sequence $b_n^*$ and continuity of $S$. By the reproducing formula $b_n^*h = \langle h, h_n \rangle_{\mathbb{H}}$ for every $h \in \mathbb{H}$, whence $\sum_n (b_n^*h)^2 < \infty$. Conversely, if $\sum_n (b_n^*b)^2 <$ $\infty$, then $h := \sum_n (b_n^* b) h_n$ is a well-defined element of $\mathbb{H}$, with $b_m^* h = b_m^* b$ for every $m$ because $b_m^* h_n = \langle h_n, h_m \rangle_{\mathbb{H}} = \delta_{mn}$. Because the linear span of the sequence $(b_n^*)$ is dense in $\mathbb{B}^*$ it follows that $b^* h = b^* b$ for every $b^*$ and hence $b = h$, which is contained in $\mathbb{H}$.

The map $\phi: \mathbb{B} \to \mathcal{R}^\infty$ defined by $b \mapsto (b_n^* b)$ is well defined and measurable. It maps $W$ onto a sequence $(Z_n) = \phi(W)$ of standard normal variables and maps $W + b$ onto the sequence $(Z_n + b_n^* b)$ of independent shifted normal variables. By Kakutani's dichotomy the latter two laws are orthogonal if $\sum_n (b_n^* b)^2 = \infty$. This implies the orthogonality of the laws of $W$ and $W + b$.

Finally we prove the existence of $(b_n^*)$ as claimed. Starting with an arbitrary dense sequence $(b_n^*)$ in $\mathbb{B}^*$, we can make this linearly independent by removing from left to right in the sequence $b_1^*, b_2^*, \ldots$ every $b_n^*$ that can be written as a linear combination of the preceding (left-over) $b_j^*$. This procedure yields a linearly independent sequence $(b_n^*)$ whose span is dense in $\mathbb{B}^*$. The random variables $b_n^* W$ are automatically linearly independent in $L_2(\Omega, \mathcal{U}, P)$, because $\sum_n \lambda_n b_n^* W = 0$, almost surely for a sequence $\lambda_n$ with finitely many nonzero elements. This implies that $\sum_n \lambda_n b_n^*$ is zero on a set with probability one under the law of $W$, and hence by continuity also on the support of this law, which is $\mathbb{B}$ by assumption. Thus we can apply the Gramm–Schmidt procedure to turn the sequence $b_n^* W$ into a sequence of standard normal variables $(Z_n)$. Then $Z_n = \sum_{i=1}^n \lambda_{i,n} b_i^* W$ for every $n$ for a triangular array of coefficients $(\lambda_{i,n})$ with $\lambda_{n,n} \neq 0$ for every $n$. The sequence $\sum_{i=1}^n \lambda_{i,n} b_i^*$ has the desired properties. □

## 4. Series representation

Suppose that the covariance kernel $K$ of the Gaussian process $W = (W_t : t \in T)$, defined on the probability space $(\Omega, \mathcal{U}, P)$, can be written in the form

$$(4.1) \qquad K(s,t) = \sum_{j=1}^\infty \lambda_j \phi_j(s) \phi_j(t)$$

for positive numbers $\lambda_1, \lambda_2, \ldots$ and arbitrary functions $\phi_j : T \to \mathcal{R}$, where the series is assumed to converge pointwise on $T \times T$. The convergence on the diagonal implies that $\sum_j \lambda_j \phi_j^2(t) < \infty$ for all $t \in T$. Then by the Cauchy–Schwarz inequality the series $\sum_{j=1}^\infty w_j \phi_j(t)$ converges absolutely for every sequence $(w_j)$ of numbers with $\sum_j w_j^2 / \lambda_j < \infty$, for every $t$, and hence defines a function from $T$ to $\mathcal{R}$. We assume that the functions $\phi_j$ are linearly independent in the sense that $\sum_j w_j \phi_j(t) = 0$ for every $t \in T$ for some sequence $(w_j)$ with $\sum_j w_j^2 / \lambda_j < \infty$ implying that $w_j = 0$ for every $j \in \mathbb{N}$.

**Theorem 4.1.** *If the covariance function $K$ of the mean-zero Gaussian process $W = (W_t : t \in T)$ can be represented as in (4.1) for numbers $\lambda_j$ and functions $\phi_j : T \to \mathcal{R}$ which satisfy $\sum_{j=1}^\infty \lambda_j \phi_j^2(t) < \infty$ for every $t \in T$ and are linearly independent as indicated, then the RKHS of the stochastic process $W$ is the set of all functions $\sum_{j=1}^\infty w_j \phi_j(t)$ with $\sum_{j=1}^\infty w_j^2 / \lambda_j < \infty$, and the inner product is given by*

$$(4.2) \qquad \Big\langle \sum_{i=1}^\infty v_i \phi_i, \sum_{j=1}^\infty w_j \phi_j \Big\rangle_{\mathbb{H}} = \sum_{j=1}^\infty \frac{v_j w_j}{\lambda_j}.$$

*Proof.* Under the condition that $\sum_{k=1}^{\infty} \lambda_k \phi_k^2(t) < \infty$ for every $t \in T$, the infinite sum defining $K(s,t)$ converges for every $(s,t) \in T \times T$, by the Cauchy–Schwarz inequality, and hence the kernel is well defined. Let $H$ be the set of all series $\sum_{k=1}^{\infty} f_k \phi_k$ when $(f_k)$ ranges over the sequences with $\sum_{k=1}^{\infty} f_k^2/\lambda_k < \infty$. (These series were noted to converge pointwise absolutely before the statement of the theorem.) By the assumed linear independence of the functions $\phi_j$, the coefficients $(f_j)$ are identifiable from the corresponding functions $\sum_j f_j \phi_j \in H$. Therefore we can define a bijection $i: H \to \ell_2$ by $i: \sum_k f_k \phi_k \mapsto (f_k/\sqrt{\lambda_k})$. The set $H$ becomes a Hilbert space under the inner product induced from $\ell_2$, which is given on the right side of (4.2), and which we denote by $\langle \cdot, \cdot \rangle_H$. We must prove that this inner product agrees with the inner product of $\mathbb{H}$ and that $H$ and $\mathbb{H}$ are the same as sets.

The function $K(s, \cdot)$ has a representation $\sum_{k=1}^{\infty} f_k \phi_k$ for $f_k = \lambda_k \phi_k(s)$, and hence is contained in $H$. It also follows that

$$\langle K(s, \cdot), K(t, \cdot) \rangle_H = \sum_{k=1}^{\infty} \frac{\lambda_k \phi_k(s) \lambda_k \phi_k(t)}{\lambda_k} = K(s,t) = \langle K(s, \cdot), K(t, \cdot) \rangle_{\mathbb{H}},$$

where the second equality follows from the series representation of $K$, and the third is (2.2). Thus the inner products of $H$ and $\mathbb{H}$ agree. We conclude that $H$ contains $\mathbb{H}$ isometrically.

The space $H$ has the reproducing property: $\langle f, K(t, \cdot) \rangle_H = f(t)$ for every $t \in T$ and $f \in H$. This follows from

$$\langle f, K(t, \cdot) \rangle_H = \left\langle \sum_k f_k \phi_k, \sum_k \lambda_k \phi_k \right\rangle_H = \sum_k \frac{f_k \lambda_k \phi_k(t)}{\lambda_k} = f(t).$$

If $f \in H$ with $f \perp \mathbb{H}$, then in particular $f \perp K(t, \cdot)$ for every $t \in T$ and hence $f(t) = 0$ by the reproducing formula. Thus $H = \mathbb{H}$. □

Series expansions of the type (4.1) are not unique, and some may be more useful than others. They may arise as an eigenvalue expansion of the operator corresponding to the covariance function. However, this is not a requirement of the proposition, which applies to arbitrary functions $\phi_j$.

**Example 4.1.** Suppose that $(T, \Theta, \nu)$ is a measurable space and

$$\iint K^2(s,t) \, d\nu(s) \, d\nu(t) < \infty.$$

Then the integral operator $K: L_2(T, \Theta, \nu) \to L_2(T, \Theta, \nu)$ defined by

$$Kf(t) = \int f(s) \, K(s,t) \, d\nu(t)$$

is compact and positive self-adjoint. Thus there exists a sequence of eigenvalues $\lambda_k \downarrow 0$ and an orthonormal system of eigenfunctions $\phi_k \in L_2(T, \Theta, \nu)$ (thus $K\phi_k = \lambda_k \phi_k$ for every $k \in \mathbb{N}$) such that (4.1) holds, where the series converges in $L_2(T \times T, \Theta \times \Theta, \nu \times \nu)$. The series $\sum_k f_k \phi_k$ now converges in $L_2(T, \Theta, \nu)$ for any sequence $(f_k)$ in $\ell_2$. By the orthonormality of the functions $\phi_k$, they are certainly linearly independent.

If the series (4.1) also converges pointwise on $T \times T$, then in particular $K(t,t) = \sum_k \lambda_k \phi_k^2(t) < \infty$ for all $t \in T$ and Theorem 4.1 shows that the RKHS is the set of all functions $\sum_k f_k \phi_k$ for sequences $(f_k)$ such that $(f_k/\sqrt{\lambda_k}) \in \ell_2$.

If the kernel is suitably regular, then we can apply the preceding with many choices of measure $\nu$, leading to different eigenfunction expansions.

If the process itself can be expanded as a series

$$W = \sum_{j=1}^{\infty} \mu_j Z_j \phi_j,$$

for a sequence of i.i.d. standard normal variables $(Z_j)$ and suitable functions $\phi_j$, where the series converges in $L_2(\Omega, \mathcal{U}, P)$, then (4.1) holds with $\lambda_j = \mu_j^2$ and the stochastic process RKHS takes the form given by the preceding proposition. The following proposition gives a Banach space version of this result.

**Theorem 4.2.** *Let $(h_i)$ be a sequence of elements in a separable Banach space $\mathbb{B}$ such that $\sum_{i=1}^{\infty} w_i h_i = 0$ for a sequence $w \in \ell_2$, where the convergence is in $\mathbb{B}$, implying that $w = 0$. Let $(Z_i)$ be an i.i.d. sequence of standard normal variables and assume that the series $W = \sum_{i=1}^{\infty} Z_i h_i$ converges almost surely in $\mathbb{B}$. Then the RKHS of $W$ as a map in $\mathbb{B}$ is given by $\mathbb{H} = \{\sum_{i=1}^{\infty} w_i h_i \colon w \in \ell_2 < \infty\}$ with squared norm $\|\sum w_i h_i\|_{\mathbb{H}}^2 = \sum_i w_i^2$.*

*Proof.* The almost sure convergence of the series $W = \sum_{i=1}^{\infty} Z_i h_i$ in $\mathbb{B}$ implies the almost sure convergence of the series $b^* W = \sum_{i=1}^{\infty} Z_i b^* h_i$ in $\mathcal{R}$, for any $b^* \in \mathbb{B}^*$. Because the partial sums of the last series are zero-mean Gaussian, the series converges also in $L_2(\Omega, \mathcal{U}, P)$. Hence for any $b^*, \underline{b}^* \in \mathbb{B}^*$,

$$\mathrm{E} b^* W \underline{b}^* W = \mathrm{E} \sum_{i=1}^{\infty} Z_i b^* h_i \sum_{i=1}^{\infty} Z_i \underline{b}^* h_i = \sum_{i=1}^{\infty} b^* h_i \underline{b}^* h_i.$$

In particular, the sequence $(b^* h_i)$ is contained in $\ell_2$ for every $b^* \in \mathbb{B}^*$, with square norm $\mathrm{E}(b^* W)^2$.

For $w \in \ell_2$ and natural numbers $m < n$, by the Hahn–Banach theorem and the Cauchy–Schwarz inequality,

$$\Big\| \sum_{m < i \le n} w_i h_i \Big\|^2 = \sup_{\|b^*\| \le 1} \Big\| \sum_{m < i \le n} w_i b^* h_i \Big\|^2$$
$$\le \sum_{m < i \le n} w_i^2 \sup_{\|b^*\| \le 1} \sum_{m < i \le n} (b^* h_i)^2.$$

As $m, n \to \infty$ the first factor on the far right tends to zero, since $w \in \ell_2$. By the first paragraph the second factor is bounded by $\sup_{\|b^*\| \le 1} \mathrm{E}(b^* W)^2 \le \mathrm{E}\|W\|^2$. Hence the partial sums of the series $\sum_i w_i h_i$ form a Cauchy sequence in $\mathbb{B}$, whence the infinite series converges.

Because $\sum_i (b^* h_i)^2$ was seen to converge, it follows that $\sum_i (b^* h_i) h_i$ converges in $\mathbb{B}$, and hence $\underline{b}^*(\sum_i (b^* h_i) h_i) = \sum_i b^* h_i \underline{b}^* h_i = \mathrm{E} b^* W \underline{b}^* W$, for any $\underline{b}^* \in \mathbb{B}^*$. This shows that $Sb^* = \sum_i (b^* h_i) h_i$ and the RKHS is not bigger than the space, as claimed.

The space would be smaller than claimed if there existed $w \in \ell_2$ that is not in the closure of the linear span of the elements $(b^* h_i)$ of $\ell_2$ when $b^*$ ranges over $\mathbb{B}^*$. We can take this $w$ without loss of generality as orthogonal to the latter collection, i.e. $\sum_i w_i b^* h_i = 0$ for every $b^* \in \mathbb{B}^*$. This is equivalent to $\sum_i w_i h_i = 0$, which has been excluded for any $w \ne 0$. □

It should be noted that the sequence $(h_i)$ in the preceding lemma consists of arbitrary elements of the Banach space, only restricted by the linear independence condition that $\sum_i w_i h_i = 0$ for $w \in \ell_2$, implying that $w = 0$ (and the convergence of the random sequence $\sum_i Z_i h_i$). Combined with an i.i.d. standard normal sequence as coefficients, this sequence turns into an *orthonormal* basis of the RKHS.

From the proof it can be seen that the linear independence is necessary. If it fails, then the RKHS is the set of linear combinations $\sum_i w_i h_i$ with $w$ restricted to the closure in $\ell_2$ of the set of sequences $(b^* h_i)$ when $b^*$ ranges over $\mathbb{B}^*$ and square norm $\sum_i w_i^2$. (Taking these linear combinations for all $w \in \ell_2$ gives the same set, but the $\ell_2$-norm should be computed for a projected $w$.)

**Example 4.2.** For $Z_0, \ldots, Z_k$ i.i.d. standard normal variables consider the polynomial process $t \mapsto \sum_{i=0}^{k} Z_i t^i/i!$ viewed as a map in (for instance) $C[0,1]$. The RKHS of this process is equal to the set of $k$th degree polynomials $P_a(t) = \sum_{i=0}^{k} a_i t^i/i!$ with square norm $\|P_a\|_\mathbb{H}^2 = \sum_{i=0}^{k} a_i^2$, i.e., the $k$th degree polynomials $P$ with square norm $\|P\|_\mathbb{H}^2 = \sum_{i=0}^{k} P^{(i)}(0)^2$.

Conversely, any Gaussian random element $W$ in a separable Banach space can be expanded in a series $W = \sum_{j=1}^{\infty} Z_j h_j$ for i.i.d. standard normal variables $Z_i$ and any orthonormal basis $(h_i)$ of its RKHS, where the series converges in the norm of the Banach space. Because we can rewrite this expansion as $W = \sum_j \|h_j\| Z_j \tilde{h}_j$, where $\tilde{h}_j = h_j/\|h_j\|$ is a sequence of norm one, the corresponding "eigenvalues" $\lambda_i$ are in this case the square norms $\|h_i\|^2$. To prove this result, recall the isometry $U: \mathbb{H} \to L_2(\Omega, \mathcal{U}, P)$ defined in (3.3).

**Theorem 4.3.** *Let $(h_i)$ be a complete orthonormal system in the RKHS $\mathbb{H}$ of a Borel measurable, zero-mean Gaussian random element $W$ in a separable Banach space $\mathbb{B}$. Then $Uh_1, Uh_2, \ldots$ is an i.i.d. sequence of standard normal variables and $W = \sum_{i=1}^{\infty} (Uh_i) h_i$, where the series converges in the norm of $\mathbb{B}$, almost surely.*

*Proof.* It is immediate from the definitions of $U$ and the RKHS that $U: \mathbb{H} \to L_2(\Omega, U, P)$ is an isometry. Because $U$ maps the subspace $S\mathbb{B}^* \subset \mathbb{H}$ into the Gaussian process $b^*W$, it maps the completion $\mathbb{H}$ of $S\mathbb{B}^*$ into the completion of the linear span of this process in $L_2(\Omega, \mathcal{U}, P)$, which consists of normally distributed variables. Because $U$ retains inner products, it follows that $Uh_1, Uh_2, \ldots$ is a sequence of i.i.d. standard normal variables.

By the definition of $U$ and its continuity, for any $b^* \in \mathbb{B}^*$,

$$b^*W = U(Sb^*) = U\Big(\sum_{i=1}^{\infty} \langle Sb^*, h_i \rangle_\mathbb{H} h_i\Big) = \sum_{i=1}^{\infty} \langle Sb^*, h_i \rangle_\mathbb{H} Uh_i = \sum_{i=1}^{\infty} b^*(h_i) Uh_i,$$

where the last equality follows from the reproducing formula (2.6) and the series converges in $L_2(\Omega, \mathcal{U}, P)$. In other words, for any $b^* \in \mathbb{B}^*$, $b^*\big(\sum_{i=1}^{n} h_i(Uh_i)\big) \equiv \sum_{i=1}^{n} (b^*h_i) Uh_i$ converges in $L_2(\Omega, \mathcal{U}, P)$ to $b^*W$. We wish to strengthen this to convergence almost surely of $W_n := \sum_{i=1}^{n} h_i(Uh_i)$ to $W$ in $\mathbb{B}$. This is an immediate consequence of the Lévy–Ito–Nisio theorem, as given in, e.g., ([9], Theorem 2.4), according to which convergence in distribution of all "marginals" $b^* \sum_{i=1}^{n} X_i$ to the marginals $b^*W$ of some Borel measurable map $W$ in a separable Banach space, for $b^* \in \mathbb{B}^*$, implies the almost sure convergence of the series $\sum_i X_i$.

An alternative proof based on a martingale argument is given in ([9], Proposition 3.6). Let $Z_1, Z_2, \ldots$ be an orthonormal basis of the closed linear span of the variables $b^*W$ in $L_2(\Omega, \mathcal{U}, P)$. Then it can be seen that, for every $n$, $\mathrm{E}(W \mid Z_1, \ldots, Z_n) = \sum_{i=1}^{n} Z_i h_i$ in a Banach space sense, for $h_i = \mathrm{E} Z_i W$. Convergence of the infinite series follows by a martingale convergence theorem for Banach space valued variables. $\square$

## 5. Support and concentration

The RKHS of a zero-mean Gaussian random element $W$ in a separable Banach space $\mathbb{B}$ is essential for an understanding of the spread of its distribution.

To begin with, the *support* of $W$, the smallest closed set $\mathbb{B}_0$ in $\mathbb{B}$ with $\mathrm{P}(W \in \mathbb{B}_0) = 1$, is the closure of the RKHS.

**Lemma 5.1.** *The support of a mean-zero Gaussian random element $W$ in a separable Banach space $\mathbb{B}$ is the closure of its RKHS in $\mathbb{B}$. It is also the closure of the set $S\mathbb{B}^*$ in $\mathbb{B}$.*

*Proof.* We first show that the probability $\mathrm{P}(\|W\| < \varepsilon)$ of an arbitrary open ball centered around 0 is positive. Let $V$ be an independent copy of $W$. Because we can cover $\mathbb{B}$ with countably many balls of radius $\varepsilon$, there exists some ball $B(h, \varepsilon)$ with positive measure under the law of $W$. The difference $B(h, \varepsilon) - B(h, \varepsilon)$ is contained in the ball of radius $2\varepsilon$ around 0. It follows that

$$\mathrm{P}\bigl(V - W \in B(0, 2\varepsilon)\bigr) \geq \mathrm{P}\bigl(V \in B(h, \varepsilon)\bigr)\mathrm{P}\bigl(W \in B(h, \varepsilon)\bigr) > 0.$$

Now $(V - W)/\sqrt{2}$ is a zero-mean Gaussian process with the same covariance function as $W$, and hence has the same distribution as $W$. It follows that $\mathrm{P}\bigl(W \in B(0, \sqrt{2}\varepsilon)\bigr) > 0$ for every $\varepsilon > 0$.

Since the distribution of $W - h$ is equivalent to the distribution of $W$ for any $h$ in the RKHS, by Lemma 3.2, it follows that $\mathrm{P}\bigl(\|W - h\| < \varepsilon\bigr) > 0$ for any $\varepsilon > 0$ and $h \in \mathbb{H}$.

This remains true for an element $h \in \mathbb{B}$ that can be approximated arbitrarily closely by elements from the RKHS. Thus the support of $W$ contains the closure of the RKHS in $\mathbb{B}$.

By the Hahn–Banach theorem this closure $\bar{\mathbb{H}}$ can be written as

$$\bar{\mathbb{H}} = \bigcap_{b^* \in \mathbb{B}^* : b^*\mathbb{H} = 0} N(b^*),$$

where $b^*\mathbb{H} = 0$ means $b^*h = 0$ for all $h \in \mathbb{H}$ and $N(b^*) = \{b \in \mathbb{B} : b^*(b) = 0\}$ is the kernel of $b^*$. If $b^*\mathbb{H} = 0$, then in particular $b^*(Sb^*) = \mathrm{E}(b^*(W))^2 = 0$, and hence $b^*(W) = 0$ almost surely. It follows that $\mathrm{P}\bigl(W \in N(b^*)\bigr) = 1$ for every $b^*$ in the display. By the preceding display the complement $\mathbb{B} - \bar{\mathbb{H}}$ is a union of the open sets $N(b^*)^c$. Because an open set in a separable metric space is Lindelöf ([6], section 10) this union can be written as a union of countably many of the sets $N(b^*)^c$. Equivalently, the intersection in the preceding display can be restricted to a suitable countable subset. It follows that $\mathrm{P}(W \in \bar{\mathbb{H}}) = 1$.

The second assertion follows, because the RKHS-norm is stronger than the norm of the containing Banach space. Completing the set $S\mathbb{B}^*$ for the RKHS-norm before taking the closure in $\mathbb{B}$ does therefore not give a bigger set. □

An inequality of [1] gives further insight in the concentration of the distribution of $W$. Let $\mathbb{H}_1$ and $\mathbb{B}_1$ be the unit balls of the RKHS and the space $\mathbb{B}$, respectively. The inequality involves the (centered) *small ball probability*

$$e^{-\phi_0(\varepsilon)} = \mathrm{P}(W \in \varepsilon\mathbb{B}_1).$$

**Theorem 5.1.** *(Borell's inequality.) For any $\varepsilon > 0$ and $M \geq 0$,*

$$\mathrm{P}\bigl(W \in \varepsilon\mathbb{B}_1 + M\mathbb{H}_1\bigr) \geq \Phi\bigl(\Phi^{-1}(e^{-\phi_0(\varepsilon)}) + M\bigr).$$

Here $\Phi$ is the cumulative distribution function of the standard normal distribution. For fixed $\varepsilon > 0$ the right side decreases as $M \to \infty$ according to the tails of the standard normal distribution. This shows that the "geometry of the concentration" of $W$ is given by the unit ball of the RKHS. Summing the small ball $\varepsilon \mathbb{B}_1$ to the multiple $M\mathbb{H}_1$ can be seen as enlarging the latter set with an $\varepsilon$-neighbourhood. In general this is necessary to capture the mass of the $W$, because the support of $W$ is the closure of the RKHS; the RKHS itself may have probability zero. For $M \to \infty$ we obtain the equality $\mathrm{P}(W \in \varepsilon\mathbb{B}_1 + \mathbb{H}) = 1$, for any $\varepsilon > 0$, which (again) shows that $W$ is supported within the closure of $\mathbb{H}$.

**Example 5.1.** For a mean-zero normal vector $W$ in $\mathbb{B} = \mathcal{R}^k$ with covariance matrix $\Sigma$, the RKHS is the range of the covariance matrix equipped with the inner product $\langle \Sigma g, \Sigma h \rangle_{\mathbb{H}} = g^T \Sigma h$. This follows, because $\mathbb{B}^* = \mathcal{R}^k$ and, for the element $g \in \mathbb{B}^*$ given by $h \mapsto h^T g$, we have $Sg = \mathrm{E}WW^T g = \Sigma g$. The inner product of the RKHS is $\langle Sg, Sh \rangle_{\mathbb{H}} = \mathrm{E}g^T W h^T W = g^T \Sigma h$.

The unit ball $\mathbb{H}_1$ is the set $\{\Sigma h : h^T \Sigma h \le 1\}$. For nonsingular $\Sigma$ this set is the ellipsoid determined by the inverse matrix $\Sigma^{-1}$, i.e., the ellipsoid determined by the level sets of the density. For singular $\Sigma$ the distribution is concentrated on a lower-dimensional subspace, and we have a similar interpretation after projection on this subspace.

Borell's inequality is often quoted as only an exponential inequality on the norm $\|W\|$, but this is in fact a consequence. The distribution of the norm $\|W\|$ of a non-zero Borel measurable Gaussian map $W$ does not have atoms (cf., [2]) and therefore has a unique median $M(W)$.

**Corollary 5.1.** *For any $x > 0$,*
$$\mathrm{P}\bigl(\|W\| - M(W) > x\bigr) \le 1 - \Phi\bigl(x/\sigma(W)\bigr).$$

*Proof.* For $\varepsilon = M(W)$ we have $\mathrm{P}(W \in \varepsilon\mathbb{B}_1) = \mathrm{P}(\|W\| \le M(W)) = 1/2$. Hence the choices $\varepsilon = M(W)$ and $M = x/\sigma(W)$ in Borell's inequality yield the inequality $\mathrm{P}(W \in M(W)\mathbb{B}_1 + (x/\sigma(W))\mathbb{H}_1) \ge \Phi(x/\sigma(W))$. Because $\mathbb{H}_1 \subset \sigma(W)\mathbb{B}_1$ by (2.5), the left side is smaller than $\mathrm{P}(W \in (M(W) + x)\mathbb{B}_1)$, which is 1 minus the left side of the corollary. □

According to Anderson's lemma (e.g., [9], p. 73, [16], p. 72, or [17], 3.11.4) a ball of fixed radius receives maximum mass of a zero-mean Gaussian distribution if centered at the origin. The following lemma gives a lower bound on the decrease in mass if the ball is centered at an element of the RKHS. The lemma is implicit in the proof of the main result in [7], and appears explicitly as (4.16) in [8].

**Lemma 5.2.** *If $h \in \mathbb{H}$, then for every Borel measurable set $C \subset \mathbb{B}$ with $C = -C$,*
$$\mathrm{P}(W - h \in C) \ge e^{-\frac{1}{2}\|h\|_{\mathbb{H}}^2} \mathrm{P}(W \in C).$$

*Proof.* By symmetry $W$ and $-W$ are identically distributed and hence $\mathrm{P}(W + h \in C) = \mathrm{P}(-W + h \in -C) = \mathrm{P}(W - h \in C)$. By Lemma 3.2,
$$\mathrm{P}(W + h \in C) = \mathrm{E}1_C(W + h) = \mathrm{E}e^{Uh - \frac{1}{2}\|h\|_{\mathbb{H}}^2} 1_C(W).$$

This is true with $-h$ instead of $h$ as well. Combining these facts yields that
$$\mathrm{P}(W - h \in C) = \tfrac{1}{2}\mathrm{E}e^{Uh - \frac{1}{2}\|h\|_{\mathbb{H}}^2} 1_C(W) + \tfrac{1}{2}\mathrm{E}e^{U(-h) - \frac{1}{2}\|-h\|_{\mathbb{H}}^2} 1_C(W)$$

$$= e^{-\frac{1}{2}\|h\|_{\mathbb{H}}^2} \mathrm{E} \cosh(Uh) 1_C(W) \geq e^{-\frac{1}{2}\|h\|_{\mathbb{H}}^2} \mathrm{P}(W \in C),$$

since $\cosh x = (e^x + e^{-x})/2 \geq 1$ for every $x$. □

The lemma with $C$ equal to the ball of radius $\varepsilon$ around 0 refers to the *noncentered small ball probabilities* $\mathrm{P}(\|W - w\| < \varepsilon)$, for every $w$ in the RKHS. Up to constants these can be completely characterized through the corresponding centered small ball probabilities and approximation of the center $w$ from the RKHS. Define

$$(5.1) \qquad \phi_w(\varepsilon) = \inf_{h \in \mathbb{H}: \|h - w\| \leq \varepsilon} \tfrac{1}{2}\|h\|_{\mathbb{H}}^2 - \log \mathrm{P}(\|W\| < \varepsilon).$$

For $w = 0$ this agrees with the negative exponent $\phi_0(\varepsilon)$ of the small ball probability $\mathrm{P}(\|W\| < \varepsilon) = e^{-\phi_0(\varepsilon)}$ defined previously. Up to constants this quantity gives the exponent of the small ball probability at center $w$.

**Lemma 5.3.** *For any $w$ in the support of $W$ and every $\varepsilon > 0$,*

$$\phi_w(\varepsilon) \leq -\log \mathrm{P}(\|W - w\| < \varepsilon) \leq \phi_w(\varepsilon/2).$$

*Proof.* For any $h \in \mathbb{H}$ with $\|h - w\| \leq \varepsilon$ we have $\|W - w\| \leq \varepsilon + \|W - h\|$ and hence $\mathrm{P}(\|W - w\| < 2\varepsilon) \geq \mathrm{P}(\|W - h\| < \varepsilon)$. The latter probability can be bounded below by $\exp(-\tfrac{1}{2}\|h\|_{\mathbb{H}}^2) \mathrm{P}(\|W\| < \varepsilon)$, in view of the preceding lemma. We conclude by optimizing over $h \in \mathbb{H}$.

The set $B_\varepsilon = \{h \in \mathbb{H}: \|h - w\| \leq \varepsilon\}$ is convex and closed in $\mathbb{H}$, because the RKHS topology is stronger than the norm topology. Therefore the (convex) map $h \mapsto \|h\|_{\mathbb{H}}^2$ attains a minimum on $B_\varepsilon$ at some point $h_\varepsilon$. Because $(1 - \lambda)h_\varepsilon + \lambda h \in B_\varepsilon$ for every $h \in B_\varepsilon$ and $0 \leq \lambda \leq 1$, it follows that $\|(1 - \lambda)h_\varepsilon + \lambda h\|_{\mathbb{H}}^2 \geq \|h_\varepsilon\|_{\mathbb{H}}^2$, which implies that $2\lambda \langle h - h_\varepsilon, h_\varepsilon \rangle_{\mathbb{H}} + \lambda^2 \|h - h_\varepsilon\|_{\mathbb{H}}^2 \geq 0$. The fact that this is true for every $0 \leq \lambda \leq 1$ can be seen to imply that $\langle h, h_\varepsilon \rangle_{\mathbb{H}} \geq \|h_\varepsilon\|_{\mathbb{H}}^2$ for every $h \in B_\varepsilon$.

By Theorem 4.3 the process $W$ can be written as $W = \sum_{i=1}^\infty (Uh_i) h_i$, for any given complete orthonormal system $h_1, h_2, \ldots$ in $\mathbb{H}$, where the series converges almost surely in norm. The truncated series $W^m = \sum_{i=1}^m (Uh_i) h_i$ takes its values in $\mathbb{H}$. If $\|W - g - w\| < \varepsilon$ and some arbitrary $g \in \mathbb{H}$, then $\|W^m - g - w\| < \varepsilon$ for sufficiently large $m$, almost surely. Equivalently, $W^m - g \in B_\varepsilon$ and hence the preceding paragraph implies that $\langle W^m - g, h_\varepsilon \rangle_{\mathbb{H}} \geq \|h_\varepsilon\|_{\mathbb{H}}^2$, eventually as $m \to \infty$, almost surely. Here $\langle W^m, h_\varepsilon \rangle_{\mathbb{H}} = \sum_{i=1}^m (Uh_i)\langle h_i, h_\varepsilon \rangle_{\mathbb{H}} = U \sum_{i=1}^m h_i \langle h_i, h_\varepsilon \rangle_{\mathbb{H}}$. By the continuity of $U$ the right side converges in $L_2(\Omega, \mathcal{U}, \mathrm{P})$ to $Uh_\varepsilon$ as $m \to \infty$, and hence almost surely along a subsequence. We conclude that $Uh_\varepsilon - \langle g, h_\varepsilon \rangle_{\mathbb{H}} \geq \|h_\varepsilon\|_{\mathbb{H}}^2$ almost surely on the event $\{\|W - g - w\| < \varepsilon\}$. In particular the choice $g = -h_\varepsilon$ yields that $Uh_\varepsilon \geq 0$ almost surely on the event $\{\|W + h_\varepsilon - w\| < \varepsilon\}$.

By Lemma 3.2,

$$\begin{aligned} \mathrm{P}(W \in w + \varepsilon \mathbb{B}_1) &= \mathrm{P}(W - h_\varepsilon \in w - h_\varepsilon + \varepsilon \mathbb{B}_1) \\ &= \mathrm{E} e^{-Uh_\varepsilon - \frac{1}{2}\|h_\varepsilon\|_{\mathbb{H}}^2} 1_{W \in w - h_\varepsilon + \varepsilon \mathbb{B}_1} \leq e^{-\frac{1}{2}\|h_\varepsilon\|_{\mathbb{H}}^2} \mathrm{E} 1_{W \in w - h_\varepsilon + \varepsilon \mathbb{B}_1}, \end{aligned}$$

by the preceding paragraph. The probability on the right side is smaller than $\mathrm{P}(W \in \varepsilon \mathbb{B}_1)$ by Anderson's lemma. □

## 6. Small ball probability and entropy

The unit ball of the RKHS not only expresses the shape of the Gaussian measure, but also allows a quantitative estimate of the small ball probability $e^{-\phi_0(\varepsilon)} = \mathrm{P}(\|W\| < \varepsilon)$ through its entropy within the Banach space.

Let $N(\varepsilon, \mathbb{H}_1, \|\cdot\|)$ be the smallest number of balls of radius $\varepsilon > 0$ needed to cover the unit ball $\mathbb{H}_1$ of the RKHS. This is bounded by the maximal number $D(\varepsilon)$ of points $h_i$ in $\mathbb{H}_1$ with $\|h_i - h_j\| \geq \varepsilon$ for $i \neq j$. Because each ball of radius $\varepsilon/2$ around a point $h_i$ has probability at least $e^{-1/2}\mathrm{P}(\|W\| < \varepsilon/2)$ by Lemma 5.2 and these balls are disjoint, it follows that $1 \geq D(\varepsilon)e^{-1/2}\mathrm{P}(\|W\| < \varepsilon/2)$, whence $D(\varepsilon)$ is finite for every $\varepsilon > 0$. This shows that the RKHS unit ball $\mathbb{H}_1$ is precompact in $\mathbb{B}$.

The following results, which were proved by [7] and [10], refine this argument, and show roughly that for regularly behaved entropy $\varepsilon \mapsto \log N(\varepsilon, \mathbb{H}_1, \|\cdot\|)$ and small ball exponent $\varepsilon \mapsto \phi_0(\varepsilon)$, and for small $\varepsilon$,

$$\log N\Big(\frac{\varepsilon}{\sqrt{\phi_0(\varepsilon)}}, \mathbb{H}_1, \|\cdot\|\Big) \asymp \phi_0(\varepsilon).$$

However, the exact statement has several constants in it.

**Lemma 6.1.** *Let $f:(0,\infty) \to (0,\infty)$ be regularly varying at zero. Then*

(i) $\log N(\varepsilon/\sqrt{2\phi_0(\varepsilon)}, \mathbb{H}_1, \|\cdot\|) \gtrsim \phi_0(2\varepsilon)$.
(ii) *If $\phi_0(\varepsilon) \lesssim f(\varepsilon)$, then $\log N(\varepsilon/\sqrt{f(e)}, \mathbb{H}_1, \|\cdot\|) \lesssim f(\varepsilon)$.*
(iii) *If $\log N(\varepsilon, \mathbb{H}_1, \|\cdot\|) \gtrsim f(\varepsilon)$, then $\phi_0(\varepsilon) \gtrsim f(\varepsilon/\sqrt{\phi_0(\varepsilon)})$.*
(iv) *If $\log N(\varepsilon, \mathbb{H}_1, \|\cdot\|) \lesssim f(\varepsilon)$, then $\phi_0(2\varepsilon) \lesssim f(\varepsilon/\sqrt{\phi_0(\varepsilon)})$.*

**Lemma 6.2.** *For $\alpha > 0$ and $\beta \in \mathcal{R}$, as $\varepsilon \downarrow 0$, $\phi_0(\varepsilon) \asymp \varepsilon^{-\alpha}(\log 1/\varepsilon)^\beta$ if and only if $\log N(\varepsilon, \mathbb{H}_1, \|\cdot\|) \asymp \varepsilon^{-2\alpha/(2+\alpha)}(\log 1/\varepsilon)^{2\beta/(2+\alpha)}$.*

## 7. RKHS under transformation

If a Gaussian process is transformed into another Gaussian process under a one-to-one, continuous, linear map, then the RKHS is transformed in parallel.

**Lemma 7.1.** *Let $T: \mathbb{B} \to \underline{\mathbb{B}}$ be a one-to-one, continuous, linear map from a separable Banach space $\mathbb{B}$ into a Banach space $\underline{\mathbb{B}}$ and let $W$ be a Borel measurable, zero-mean Gaussian random element in $\mathbb{B}$ with RKHS $\mathbb{H}$. Then the RKHS of the Gaussian random element $TW$ in $\underline{\mathbb{B}}$ is equal to $T\mathbb{H}$ and $T: \mathbb{H} \to \underline{\mathbb{H}}$ is an isometry for the RKHS-norms.*

*Proof.* Let $T^*: \underline{\mathbb{B}}^* \to \mathbb{B}^*$ be the adjoint of $T$. The RKHS $\underline{\mathbb{H}}$ of $TW$ is by definition the completion of the set of Pettis integrals

$$\underline{S}\underline{b}^* = \mathrm{E}(TW)\underline{b}^*(TW) = T(\mathrm{E}W\underline{b}^*(TW)) = TST^*\underline{b}^*,$$

for the inner product

$$\langle \underline{S}\underline{b}_1^*, \underline{S}\underline{b}_2^* \rangle_{\underline{\mathbb{H}}} = \mathrm{E}\underline{b}_1^*(TW)\underline{b}_2^*(TW) = \mathrm{E}(T^*\underline{b}_1^*W)(T^*\underline{b}_2^*W) = \langle ST^*\underline{b}_1^*, ST^*\underline{b}_2^* \rangle_{\mathbb{H}}.$$

It follows that the element $\underline{S}\underline{b}^*$ of $\underline{\mathbb{H}}$ is the image under $T$ of the element $ST^*\underline{b}^*$ of $\mathbb{H}$, and its norm is the same: $\|\underline{S}\underline{b}^*\|_{\underline{\mathbb{H}}} = \|ST^*\underline{b}^*\|_{\mathbb{H}}$. Thus $T: ST^*\underline{\mathbb{B}}^* \subset \mathbb{H} \to \underline{\mathbb{H}}$ is an isometry for the RKHS-norms. It extends by continuity to a linear map from the completion $\mathbb{H}_0$ of $ST^*\underline{\mathbb{B}}^*$ in $\mathbb{H}$ to $\underline{\mathbb{H}}$. Because $T$ is continuous for the norm of $\mathbb{B}$, this extension agrees with $T$. Because $T: ST^*\underline{\mathbb{B}}^* \to \underline{S}\underline{\mathbb{B}}^*$ is onto, $T$ is an isometry for the RKHS-norms, and $\mathbb{H}_0$ and $\underline{\mathbb{H}}$ are by definition the completions of $ST^*\underline{\mathbb{B}}^*$ and $\underline{S}\underline{\mathbb{B}}^*$, we have that $T: \mathbb{H}_0 \to \underline{\mathbb{H}}$ is an isometry onto $\underline{\mathbb{H}}$. It remains to be shown that $\mathbb{H}_0 = \mathbb{H}$.

Because $T$ is one-to-one, the range $T^*\underline{\mathbb{B}}^*$ of its adjoint is weak-* dense in $\mathbb{B}^*$ ([12], Corollary 4.12). By Lemma 2.2 the map $S^{:}\mathbb{B}^* \to \mathbb{H}$ is continuous relative to the weak-* and RKHS topologies. Combined this yields that $S(T^*\underline{\mathbb{B}}^*)$ is dense in $S\mathbb{B}^*$ for the RKHS-norm of $\mathbb{H}$ and hence is dense in $\mathbb{H}$.

Taken together the preceding shows that $T: \mathbb{H} \to \underline{\mathbb{H}}$ is an isometry onto $\underline{\mathbb{H}}$. □

## 8. RKHS relative to different norms

A stochastic process $W$ can often be viewed as a map into several Banach spaces. For instance, a process indexed by the unit interval with continuous sample paths is a Borel measurable map in the space $C[0,1]$, but also in the space $L_2[0,1]$; a process with continuously differentiable sample paths is a map in $C[0,1]$, but also in $C^1[0,1]$. The RKHS obtained from using a weaker Banach space is typically the same.

**Lemma 8.1.** *Let $(\mathbb{B}, \|\cdot\|)$ be a separable Banach space on $\mathbb{B}$ and let $\|\cdot\|'$ be a norm on $\mathbb{B}$ with $\|b\|' \leq \|b\|$. Then the RKHS of a Borel measurable zero-mean Gaussian random element in $(\mathbb{B}, \|\cdot\|)$ is the same as the RKHS of this map viewed in the completion of $\mathbb{B}$ under $\|\cdot\|'$.*

*Proof.* Let $\mathbb{B}'$ be the completion of $\mathbb{B}$ relative to $\|\cdot\|'$. The assumptions imply that the identity map $I:(\mathbb{B}, \|\cdot\|) \to (\mathbb{B}', \|\cdot\|')$ is continuous, linear and one-to-one. The proposition therefore is a consequence of Lemma 7.1. □

**Example 8.1.** Let $W$ be a mean zero Gaussian process indexed by the unit interval $[0,1]$ with covariance function $K(s,t) = \mathrm{E} W_s W_t$.

If $W$ has continuous sample paths, then it is a random element in $C[0,1]$. The RKHS of $W$ viewed as a random element in $C[0,1]$ is the completion of the linear span of the functions $K(t, \cdot)$ under the inner product (2.2).

If $W$ is a measurable process and $\int_0^1 W_s^2 \, ds < \infty$ surely, then $W$ is a random element in $L_2[0,1]$. The dual space of $L_2[0,1]$ consists of the maps $g \mapsto \int g(s) f(s) \, ds$ for $f$ ranging over $L_2[0,1]$, and $Sf(t) = \mathrm{E} W_t \int W_s f(s) \, ds = \int K(s,t) f(s) \, ds$. Therefore, the RKHS of $W$ viewed as a random element in $L_2[0,1]$ is the completion of the linear span of the functions $t \mapsto \int K(s,t) f(s) \, ds$ for $f$ ranging over $L_2[0,1]$ under the inner product $\langle Sf, Sg \rangle_{\mathbb{H}} = \int \int K(s,t) f(s) g(t) \, ds \, dt$.

If $W$ has continuous sample paths, then its covariance kernel is continuous, and it can be shown by direct arguments that the two RKHSs agree. This also follows from the preceding lemma.

## 9. RKHS under independent sums

If a given Gaussian prior misses certain desirable "directions" in its RKHS, then these can be filled in by adding independent Gaussian components in these directions. A closed linear subspace $\mathbb{B}_0 \subset \mathbb{B}$ of a Banach space $\mathbb{B}$ is *complemented* if there exists a closed linear subspace $\mathbb{B}_1$ with $\mathbb{B} = \mathbb{B}_0 + \mathbb{B}_1$ and $\mathbb{B}_0 \cap \mathbb{B}_1 = \{0\}$.

**Lemma 9.1.** *Let $V$ and $W$ be independent Borel measurable, zero-mean, jointly Gaussian maps from a given probability space into a separable Banach space with supports $\mathbb{B}^V$ and $\mathbb{B}^W$ such that $\mathbb{B}^V \cap \mathbb{B}^W = \{0\}$ and the subspace $\mathbb{B}^V$ is complemented by a subspace that contains $\mathbb{B}^W$. Then the RKHS of $V + W$ is the direct sum of the RKHSs of $V$ and $W$ and the RKHS norms satisfy $\|h^V + h^W\|_{\mathbb{H}^{V+W}}^2 = \|h^V\|_{\mathbb{H}^V}^2 + \|h^W\|_{\mathbb{H}^W}^2$.*

*Proof.* By the independence of $V$ and $W$ the Pettis integral $S^{V+W}b^* = \mathrm{E}(W + V)b^*(V + W)$ can be written as $S^{V+W}b^* = S^V b^* + S^W b^*$. The assumptions of trivial intersection $\mathbb{B}^V \cap \mathbb{B}^W = \{0\}$ and of complementation of $\mathbb{B}^V$ entail that there exists a continuous linear map $\Pi \colon \mathbb{B} \to \mathbb{B}^V$ such that $\Pi b = b$ if $b \in \mathbb{B}^V$ and $\Pi b = 0$ if $b \in \mathbb{B}^W$; (cf., [6], 29.2). Then $b^* \circ \Pi \in \mathbb{B}^*$, $\Pi V = V$ and $\Pi W = 0$ almost surely, whence $S^W(b^* \circ \Pi) = 0$ and hence $S^{V+W}(b^* \circ \Pi) = S^V b^*$. It follows that $S^V \mathbb{B}^* \subset S^{V+W}\mathbb{B}^*$ and by symmetry $S^W \mathbb{B}^* \subset S^{V+W}\mathbb{B}^*$. Also, for any $b_1^*, b_2^* \in \mathbb{B}^*$,

$$\begin{aligned}
\langle S^V b_1^*, S^W b_2^* \rangle_{\mathbb{H}^{V+W}} &= \langle S^{V+W}b_1^* \circ \Pi, S^{V+W}b_2^* \circ (I - \Pi) \rangle_{\mathbb{H}^{V+W}} \\
&= \mathrm{E}\bigl(b_1^* \circ \Pi(V+W)\bigr)\bigl(b_2^* \circ (I - \Pi)(V+W)\bigr) \\
&= \mathrm{E}(b_1^* V)(b_2^* W) = 0.
\end{aligned}$$

We conclude that $S^V \mathbb{B}^* \perp S^W \mathbb{B}^*$ in $\mathbb{H}^{V+W}$, so that $\mathbb{H}^{V+W}$ is the direct (orthogonal) sum of $\mathbb{H}^V$ and $\mathbb{H}^W$. Furthermore $\|S^{V+W}b^*\|^2_{\mathbb{H}^{V+W}} = \mathrm{E}(b^*(V+W))^2 = \mathrm{E}(b^*V)^2 + \mathrm{E}(b^*W)^2$. □

By the Hahn–Banach theorem the assumption of complementation is certainly satisfied as soon as one of the supports of $V$ and $W$ is finite-dimensional.

The assumption that $\mathbb{B}^V \cap \mathbb{B}^W = \{0\}$ can be interpreted as requiring "linear independence" rather than some form of orthogonality of the supports of $V$ and $W$. The stochastic independence of $V$ and $W$ translates the linear independence into orthogonality in the RKHS of $V + W$.

The assumption requires trivial intersection of the supports of the variables $V$ and $W$, rather than of sets that carry probability one. Because the RKHS is independent of the norm (Lemma 8.1) the closure operation involved in computing the support may be taken for the strongest norm which is defined on the random elements.

The assumption that $\mathbb{B}^V \cap \mathbb{B}^W = \{0\}$ cannot be removed. For instance, if $V = \sum_i \mu_i Z_i \psi_i$ and $W = \sum_i \mu'_i Z'_i \psi_i$ are series expansions with independent standard normal variables $(Z_i), (Z'_i)$ on a common basis $(\psi_i)$, then the sum process can be written $V + W = \sum_i \mu''_i Z''_i \psi_i$ for $\mu''_i = \sqrt{\mu_i^2 + (\mu'_i)^2}$ and $Z''_i$ independent standard normal variables. The RKHS of $V + W$ is then the set of series $\sum_i w_i \psi_i$ with coefficients $(w_i)$ satisfying $\sum_i (w_i/\mu''_i)^2 < \infty$ (see Section 4). Thus the RKHS is not an orthogonal sum and, asymptotically as $i \to \infty$, the eigenvalues $(\mu''_i)^2$, which determine the presence of the directions $\psi_i$ in the RKHS, are determined by the slowest of the two sequences $\mu_i$ and $\mu'_i$. If $\mu_i/\mu'_i \to 0$, then the RKHS of $V + W$ is essentially the same as the RKHS of $W$.

## 10. Examples

The RKHS of standard Brownian motion, viewed as a random element in $C[0,1]$, is well known to be the set

$$(10.1) \qquad \bigl\{f \colon [0,1] \to \mathcal{R}, f \in AC, f(0) = 0, \int f'(t)^2 \, dt < \infty\bigr\},$$

where $f \in AC$ is the assumption that $f$ is absolutely continuous. The RKHS inner product is

$$\langle f, g \rangle_{\mathbb{H}} = \int_0^1 f'(t) g'(t) \, dt.$$

**Lemma 10.1.** *The RKHS of a standard Brownian motion $W$ on $[0,1]$ is given by (10.1) with the inner product as indicated.*

*Proof.* We use the definition of the RKHS in Section 2.1 and the fact that the covariance kernel of Brownian motion is given by $s \wedge t$. The RKHS is the completion of the linear span of the functions $t \mapsto s \wedge t$ as $s$ ranges over $[0, 1]$, under the inner product determined by

$$\langle s_1 \wedge \cdot, s_2 \wedge \cdot \rangle_\mathbb{H} = s_1 \wedge s_2 = \int (s_1 \wedge t)'(s_2 \wedge t)'\, dt,$$

where the prime denotes differentiation relative to $t$, in the sense of absolute continuity.

The linear span of the functions $t \mapsto s \wedge t$ contains every function that is 0 at 0, continuous, and piecewise linear on a partition $0 = s_0 < s_1 \cdots < s_N = 1$. Indeed to obtain such a function with slopes $\alpha_1, \ldots, \alpha_N$ on the intervals $(s_0, s_1), \ldots, (s_{N-1}, s_N)$, first determine the coefficient of $s_N \wedge \cdot$ to have a correct slope on $(s_{N-1}, s_N)$, next determine the coefficient of $s_{N-1} \wedge \cdot$ to have a correct slope on $(s_{N-2}, s_{N-1})$, etc. The derivatives of these functions are piecewise constant, and the set of piecewise constant functions is dense in $L_2[0, 1]$. □

Given the RKHS of Brownian motion it is now easy to derive the RKHS of several processes related to it.

- To release Brownian motion at zero, we may start it at an independent standard normal variable $Z$, giving the process $t \mapsto Z + W_t$. The RKHS of the constant process $t \mapsto Z$ are the constant functions, which have trivial intersection with the RKHS of Brownian motion. A given function $f: [0, 1] \to \mathcal{R}$ can be decomposed as $f = f(0) + (f - f(0))$, where the second part is in the RKHS of Brownian motion if it is absolutely continuous with square integrable derivative. By Lemma 9.1, the RKHS of $Z + W$ is the set of all absolutely continuous functions $f: [0, 1] \to \mathcal{R}$ equipped with the inner product $\langle f, g \rangle_\mathbb{H} = f(0)g(0) + \int f'(s)g'(s)\, ds$.
- To smooth Brownian motion we may consider its $k$-fold integral $I_{0+}^k W$, where $(I_{0+}^1 f)(t) = \int_0^t f(s)\, ds$ and $I_{0+}^k = I_{0+}^{k-1} I_{0+}^1$. Taking a primitive is a continuous, linear, one-to-one map from $C[0, 1] \to C[0, 1]$, and hence by Lemma 7.1 the RKHS of $I_{0+}^k W$ is the set of functions $I_{0+}^k f$ for $f$ in the RKHS of Brownian motion, equipped with the inner product $\langle I_{0+}^k f, I_{0+}^k g \rangle_\mathbb{H} = \int_0^1 f'(s)g'(s)\, ds$. This space can be described simply as the set of all functions $f: [0, 1] \to \mathcal{R}$ that are $k$-times differentiable with an absolutely continuous $k$th derivative with square-integrable $f^{(k+1)}$, equipped with the inner product $\langle f, g \rangle_\mathbb{H} = \int_0^1 f^{(k+1)}(s)g^{(k+1)}(s)\, ds$.
- The sample paths of $k$-fold integrated Brownian motion $I_{0+}^k W$ have $k$ vanishing derivatives at zero, which negatively affects its approximation properties to smooth functions. (See Example 10.1 below.) We can release the derivatives by adding a polynomial and considering the process $t \mapsto \sum_{i=0}^k Z_i t^i/i! + (I_{0+}^k W)_t$, for $Z_0, \ldots, Z_k$ i.i.d. standard normal variables, independent of $W$. The supports of the polynomial process $t \mapsto \sum_{i=0}^k Z_i t^i/i!$ and $I_{0+}^k W$ in $C[0, 1]$ do not have a trivial intersection, and hence we cannot apply Lemma 9.1 in that setting. However, we may consider these processes as Borel measurable random elements in the space $C^{(k)}[0, 1]$ of $k$-times differentiable functions, equipped with the norm $\|f\|_k = \|f\|_\infty + \|f^{(k)}\|_\infty$. According to Lemma 8.1, this does not change the RKHS. The support of the process $I_{0+}^k W$ in $C^{(k)}[0, 1]$ contains only functions with $k$ vanishing derivatives at 0, and hence does have trivial intersection with the support of the polynomial process $t \mapsto \sum_{i=0}^k Z_i t^i/i!$, which is the set of $k$th degree polynomials. Applied in this setting Lemma 9.1 yields that the RKHS of the process $t \mapsto \sum_{i=0}^k Z_i t^i/i! + (I_{0+}^k W)_t$ is the set of functions $f: [0,1] \to \mathcal{R}$ that are $k$-times differentiable with an absolutely continuous $k$th derivative with square-integrable $f^{(k+1)}$, equipped with the inner product $\langle f, g \rangle_{\mathbb{H}} = \sum_{i=0}^{k} f^{(i)}(0) g^{(i)}(0) + \int_0^1 f^{(k+1)}(s) g^{(k+1)}(s)\, ds$. To see the latter, note that any $f$ can be uniquely written as $f = P_k + (f - P_k)$ for $P_k(t) = \sum_{i=0}^{k} f^{(i)}(0) t^i / i!$, the $k$th degree Taylor polynomial and $f - P_k$ a function with $k$ vanishing derivatives at zero. The polynomial $P_k$ is contained in the RKHS of the polynomial process $t \mapsto \sum_{i=0}^{k} Z_i t^i / i!$ with square RKHS-norm $\sum_{i=0}^{k} P_k^{(i)}(0)^2$ by Example 4.2, and the function $f - P_k$ is contained in the RKHS of $I_{0+}^k W$ by the preceding.

The preceding can be extended to fractional integrals of Brownian motion. Rather than studying the fractional integral operator in detail, we give a direct derivation of the RKHSs. For $\alpha > 0$ and $W$ a standard Brownian motion the *Riemann–Liouville process with Hurst parameter* $\alpha > 0$ is defined as

$$R_t^\alpha = \int_0^t (t-s)^{\alpha-1/2}\, dW_s, \quad t \geq 0.$$

The process $R^\alpha$ is a centered Gaussian process with continuous sample paths. It can be viewed as a multiple of the $(\alpha + 1/2)$-fractional integral of the "derivative $dW$ of Brownian motion". For $\alpha > 0$ and a (deterministic) measurable function $f$ on $[0,1]$ the (left-sided) *Riemann–Liouville fractional integral* of $f$ of order $\alpha$ (if it exists) is defined as (cf. [13])

$$I_{0+}^\alpha f(t) = \frac{1}{\Gamma(\alpha)} \int_0^t (t-s)^{\alpha-1} f(s)\, ds.$$

For $\alpha$ a natural number, the function $I_{0+}^\alpha f$ is just the $\alpha$-fold iterated integral of $f$, and for $\alpha > 1/2$ the Rieman–Liouville process is equal to $\Gamma(\alpha + 1/2) I_{0+}^{\alpha-1/2} W$ for $I_{0+}^\alpha$ the fractional integral.

**Lemma 10.2.** *The RKHS of the Riemann–Liouville process with parameter $\alpha > 0$ viewed as a random element in $C[0,1]$ is $\mathbb{H} = I_{0+}^{\alpha+1/2}(L_2[0,1])$ and the RKHS-norm is given by*

$$\|I_{0+}^{\alpha+1/2} f\|_{\mathbb{H}} = \frac{\|f\|_2}{\Gamma(\alpha + 1/2)}.$$

*Proof.* We use the characterization of the RKHS as the completion of the functions (2.1) under the inner product (2.2). With $f_s$ the function defined by $f_s(u) = (s-u)_+^{\alpha-1/2}$, we have, for all $s, t \geq 0$ and $\langle \cdot, \cdot \rangle_2$ the inner product of $L_2[0,1]$,

$$E R_t^\alpha R_s^\alpha = \int (t-u)_+^{\alpha-1/2} (s-u)_+^{\alpha-1/2}\, du = \langle f_t, f_s \rangle_2 = \Gamma(\alpha + 1/2) I_{0+}^{\alpha+1/2} f_s(t).$$

Hence every simple element of $\mathbb{H}$ of the form (2.1) is given by $I_{0+}^{\alpha+1/2} f$ for some $f \in L_2[0,1]$. Moreover, the inner product (2.2) of two such elements $I_{0+}^{\alpha+1/2} f$ and $I_{0+}^{\alpha+1/2} g$ is given by

(10.2) $$\left\langle I_{0+}^{\alpha+1/2} f, I_{0+}^{\alpha+1/2} g \right\rangle_{\mathbb{H}} = \frac{\langle f, g \rangle_2}{\Gamma^2(\alpha + 1/2)}.$$

It follows that the RKHS $\mathbb{H}$ is a subspace of the Hilbert space obtained by endowing $I_{0+}^{\alpha+1/2}(L_2[0,1])$ with the inner product (10.2). To prove the converse inclusion, suppose that $g \in I_{0+}^{\alpha+1/2}(L_2[0,1])$ is orthogonal to $\mathbb{H}$. Then $g = I_{0+}^{\alpha+1/2} f$ for some $f \in L_2[0,1]$ and $g$ is, in particular, orthogonal to every element $I_{0+}^{\alpha+1/2} f_t$ of $\mathbb{H}$. Hence, for every $t \in [0,1]$,

$$0 = \left\langle I_{0+}^{\alpha+1/2} f, I_{0+}^{\alpha+1/2} f_t \right\rangle_{\mathbb{H}} = \frac{\langle f, f_t \rangle_2}{\Gamma^2(\alpha+1/2)} = \frac{I_{0+}^{\alpha+1/2} f(t)}{\Gamma(\alpha+1/2)}.$$

The injectivity of the operator $I_{0+}^{\alpha+1/2}: L_2[0,1] \to L_2[0,1]$ (see [13], Theorem 13.1) then implies that $f = 0$, whence $g = 0$. We conclude that $\mathbb{H} = I_{0+}^{\alpha+1/2}(L_2[0,1])$, and the inner product on $\mathbb{H}$ is given by (10.2). $\square$

**Example 10.1.** The process $t \mapsto Z + \int_0^t W_s \, ds$, for $Z$ a standard normal variable and $W$ an independent Brownian motion, has sample paths of regularity $3/2$ and can take any value at 0, but the derivative at 0 is 0. We shall show that the latter makes the process inappropriate as a prior model for $3/2$-smooth functions.

By similar arguments as before the RKHS $\mathbb{H}$ of the process can be seen to be the set of all functions $h: [0,1] \to \mathcal{R}$ with absolutely continuous derivative such that $\int_0^1 h''(s)^2 \, ds < \infty$ and $h'(0) = 0$, with square norm

$$\|h\|_{\mathbb{H}}^2 = \|h''\|_2^2 + h(0)^2.$$

We shall show that for the identity function *id* we have

$$\inf\{\|h\|_{\mathbb{H}}^2 : h \in \mathbb{H}, \|h - id\|_\infty < \varepsilon\} \gtrsim \frac{1}{\varepsilon}.$$

This may be contrasted with the approximation by the RKHS of the process $t \mapsto Z_0 + Z_1 t + \int_0^t W_s \, ds$, which is of order $(1/\varepsilon)^{2/3}$ for every function in $C^{3/2}[0,1]$ (see [15]).

To prove the claim note that $\|h - id\|_\infty < \varepsilon$ implies that $h(3\varepsilon) - h(0) > \varepsilon$. Therefore the quantity in the display is bounded below by

$$\inf\left\{\int_0^{3\varepsilon} h''(s)^2 \, ds : h(3\varepsilon) - h(0) > \varepsilon, h'(0) = 0\right\}.$$

For a given $h$ as in the display we can define $g$ by

$$g(y) = \frac{h(3\varepsilon y) - h(0)}{\varepsilon}.$$

Then $g'(y) = 3h'(3\varepsilon y)$, $g''(y) = 9h''(3\varepsilon y)\varepsilon$, and

$$g(0) = 0, \qquad g'(0) = 0, \qquad g(1) > 1.$$

$$\int_0^{3\varepsilon} h''(s)^2 \, ds = \int_0^{3\varepsilon} g''(s/(3\varepsilon))^2 \frac{1}{(9\varepsilon)^2} \, ds = \int_0^1 g''(u)^2 \frac{1}{27\varepsilon} \, du.$$

Thus the preceding display is bigger than

$$\left(\frac{1}{27\varepsilon}\right) \inf\left\{\int_0^1 g''(u)^2 \, du : g(1) > 1, g(0) = g'(0) = 0\right\}.$$

The infimum is nonzero, because $g'' = 0$ implies that $g$ is a linear function, hence identically 0 because $g(0) = g'(0) = 0$, contradicting $g(1) > 1$.